\documentstyle[11pt]{article}
\title{Some examples related to the Deligne-Simpson problem
\footnote{Research partially supported 
by INTAS grant 97-1644}}
\author{Vladimir Petrov Kostov\\ \\ \hspace{7cm}
{\sl To the memory of my mother}\\ \\Appendix by Ofer Gabber} 
\date{}
\newtheorem{tm}{Theorem}
\newtheorem{lm}[tm]{Lemma}

\newtheorem{prop}[tm]{Proposition}
\newtheorem{rem}[tm]{Remark}
\newtheorem{rems}[tm]{Remarks}
\newtheorem{defi}[tm]{Definition}
\newtheorem{ex}[tm]{Example}
\newtheorem{conjecture}[tm]{Open questions}
\begin{document}
\maketitle 

\section{Introduction}
\subsection{The Deligne-Simpson problem}

In the present paper we consider some examples relative to the 
{\em Deligne-Simpson problem (DSP)} which is formulated like this: 

{\em Give necessary and sufficient conditions upon the choice of the $p+1$ 
conjugacy classes $c_j\subset gl(n,{\bf C})$, resp. 
$C_j\subset GL(n,{\bf C})$, so that there exist irreducible $(p+1)$-tuples of 
matrices $A_j\in c_j$, $A_1+\ldots +A_{p+1}=0$, resp. of matrices 
$M_j\in C_j$, $M_1\ldots M_{p+1}=I$.}

By definition, the {\em weak DSP} is the DSP in which the requirement of 
irreducibility is replaced by the weaker requirement the centralizer of the 
$(p+1)$-tuple of matrices to be trivial. 

The matrices $A_j$, resp. $M_j$, are interpreted as {\em matrices-residua} of 
{\em Fuchsian} systems on Riemann's sphere (i.e. linear systems of ordinary 
differential equations with logarithmic poles), resp. as 
{\em monodromy operators} 
of {\em regular} systems on Riemann's sphere (i.e. linear systems of ordinary 
differential equations with moderate growth rate of the solutions at the 
poles). Fuchsian systems are a particular case of regular ones. 
By definition, the monodromy operators generate the 
{\em monodromy group} of a regular system. 

In the multiplicative version (i.e. for matrices $M_j$) 
the classes $C_j$ are interpreted as 
{\em local monodromies} around the poles and the problem admits the 
interpretation: 

{\em For what $(p+1)$-tuples of local monodromies do there exist  
monodromy groups with such local monodromies.} 

\begin{rems}
1) Suppose that $A_j$ denotes a matrix-residuum and that $M_j$ denotes the 
corresponding monodromy operator of a Fuchsian system. Then 
in the absence of non-zero integer differences 
between the eigenvalues of $A_j$ the operator $M_j$ is conjugate to 
$\exp (2\pi iA_j)$. 

2) In what follows the sum of the matrices $A_j$ is always presumed to be 
0 and the product of the matrices $M_j$ is always presumed to be $I$.
\end{rems}

\subsection{The aim of this paper}

For a conjugacy class $C$ in $GL(n,{\bf C})$ or $gl(n,{\bf C})$ denote by 
$d(C)$ its dimension (which is always even). Set $d_j:=d(c_j)$ 
(resp. $d(C_j)$). 

For fixed conjugacy classes $C_j$ consider the variety 

\[ {\cal V}=\{ (M_1,\ldots ,M_{p+1})~|~M_j\in C_j~,~M_1\ldots M_{p+1}=I\}~.\] 
This variety might contain $(p+1)$-tuples with non-trivial centralizers as 
well as with trivial ones. It might contain only the former or only the 
latter. 

\begin{prop}\label{smooth}
At a $(p+1)$-tuple with trivial centralizer the variety ${\cal V}$ is smooth 
and of dimension $d_1+\ldots +d_{p+1}-n^2+1$. 
\end{prop}

\begin{rem}
The proposition is proved at the end of the subsection. 
A similar statement is true for matrices $A_j$.
\end{rem}

For generic eigenvalues (the precise definition is given in the next 
section) the variety ${\cal V}$ contains only irreducible $(p+1)$-tuples and 
its dimension remains the same when the eigenvalues of the conjugacy classes 
are changed but not the Jordan normal forms which they define. We call its 
dimension for generic eigenvalues the {\em expected} one. 

The aim of the present paper is to consider some examples of varieties 
${\cal V}$ for non-generic eigenvalues. 
In the first and in the fifth of them (see Sections~\ref{ex1} and \ref{ex5})
dim${\cal V}$ is higher than the expected 
one. In the first example we discuss the stratified structure of ${\cal V}$ 
and we show that ${\cal V}$ consists only of $(p+1)$-tuples with non-trivial 
centralizers. The latter fact is true for the fifth example as well.  

In the second example (see Section~\ref{ex2}) 
the eigenvalues are not generic and the variety 
${\cal V}$ contains at the same time $(p+1)$-tuples with trivial and ones 
with non-trivial centralizers. The dimension of ${\cal V}$ is the 
expected one. 

In the third example (see Section~\ref{ex3}) 
the variety ${\cal V}$ contains no $(p+1)$-tuples with trivial 
centralizers but its dimension equals the expected one.

In the fourth example (see Section~\ref{ex4}) there is coexistence in 
${\cal V}$ of $(p+1)$-tuples with trivial 
centralizers and of $(p+1)$-tuples with non-trivial ones. The dimension of 
${\cal V}$ at the former (i.e. the expected dimension) is lower than the 
dimension at the latter.

In the first and third examples the closure of ${\cal V}$ 
(topological and algebraic) contains also $(p+1)$-tuples in which some of the 
matrices $M_j$ belong not to $C_j$ but to their closures, i.e. the eigenvalues 
are the necessary ones but the Jordan structure is ``less generic''.   

Similar examples exist for matrices $A_j$ as well. Before beginning with the 
examples we recall some known facts in the next section.

{\bf Proof of Proposition~\ref{smooth}:}

It suffices to prove the proposition in the case when 
$C_j\subset SL(n,{\bf C})$. 
The variety ${\cal V}$ is the intersection in 
$C_1\times \ldots \times C_p\times SL(n,{\bf C})$ of 
the graph of the mapping 

\[ C_1\times \ldots \times C_p\rightarrow SL(n,{\bf C})~,~
(M_1,\ldots ,M_p)\mapsto (M_1\ldots M_p)^{-1}\]
and of the variety ${\cal C}=C_1\times \ldots \times C_{p+1}$. To prove that 
${\cal V}$ is smooth 
it suffices to prove that the intersection is transversal, i.e. the 
sum of the tangent spaces to the graph (which is the space 
$\{ \sum _{j=1}^p[M_j,X_j],X_j\in sl(n,{\bf C})\}$) and the 
one to ${\cal C}$ (it equals 
$\{ [M_{p+1},X_{p+1}],X_{p+1}\in sl(n,{\bf C})\}$) is $sl(n,{\bf C})$. 
This follows from  

\begin{prop}\label{commut}
The $(p+1)$-tuple of matrices $R_j\in gl(n,{\bf C})$ is with trivial 
centralizer if and only if the map 
$(gl(n,{\bf C}))^{p+1}\rightarrow sl(n,{\bf C})$, 
$(X_1,\ldots ,X_{p+1})\mapsto \sum _{j=1}^{p+1}[R_j,X_j]$ is surjective.
\end{prop}

The dimension of ${\cal V}$ is the one of $C_1\times \ldots \times C_p$, i.e. 
$d_1+\ldots +d_p$, diminished by the codimension of ${\cal C}$ in 
$C_1\times \ldots \times C_p\times SL(n,{\bf C})$, i.e. by $n^2-1-d_{p+1}$. 
Hence, dim${\cal V}=d_1+\ldots +d_{p+1}-n^2+1$.$\hspace{1cm}\Box$

{\bf Proof of Proposition~\ref{commut}:}

The map is not surjective exactly if the image of every map 
$X_j\mapsto [R_j,X_j]$ belongs to one and the same linear subspace of 
$sl(n,{\bf C})$, i.e. one has tr$(D[R_j,X_j])=0$ for some matrix 
$0\neq D\in sl(n,{\bf C})$ for $j=1,\ldots ,p+1$ and identically in the 
entries of $X_j$. One has tr$(D[R_j,X_j])=$tr$([D,R_j]X_j)$ which implies that 
$[D,R_j]=0$ for all $j$ -- a contradiction with the triviality of 
the centralizer.$\hspace{1cm}\Box$

\section{Some known facts}

We expose here some facts which are given in some more detail in \cite{Ko}. 
For a matrix $Y$ from the conjugacy class $C$ in $GL(n,{\bf C})$ or 
$gl(n,{\bf C})$ set 
$r(C):=\min _{\lambda \in {\bf C}}{\rm rank}(Y-\lambda I)$. The integer 
$n-r(C)$ is the maximal number of Jordan blocks of $J(Y)$ with one and the 
same eigenvalue. Set $r_j:=r(c_j)$ 
(resp. $r(C_j)$). The quantities 
$r(C)$ and $d(C)$ depend only on the Jordan normal form of $Y$. 

\begin{defi}
A {\em Jordan normal form (JNF) of size $n$} is a family 
$J^n=\{ b_{i,l}\}$ ($i\in I_l$, $I_l=\{ 1,\ldots ,s_l\}$, $l\in L$) of 
positive integers $b_{i,l}$ 
whose sum is $n$. The index $l$ is the one of an eigenvalue and the index 
$i$ is the one of a Jordan block with the $l$-th eigenvalue; all eigenvalues 
are presumed distinct. An $n\times n$-matrix 
$Y$ has the JNF $J^n$ (notation: $J(Y)=J^n$) if to its distinct 
eigenvalues $\lambda _l$, $l\in L$, there belong Jordan blocks of sizes 
$b_{i,l}$. We usually assume that for each fixed $l$ the numbers $b_{i,l}$ 
form a non-increasing sequence.
\end{defi}

\begin{prop}\label{d_jr_j}
(C. Simpson, see \cite{Si}.) The 
following couple of inequalities is a necessary condition for the existence 
of irreducible $(p+1)$-tuples of matrices $M_j$:

\[ d_1+\ldots +d_{p+1}\geq 2n^2-2~~~~~(\alpha _n)~~,~~~~~
{\rm for~all~}j,~r_1+\ldots +\hat{r}_j+\ldots +r_{p+1}\geq n~~~~~
(\beta _n)~~~.\]
\end{prop}

\begin{rem}
The conditions are necessary for the existence of irreducible 
$(p+1)$-tuples of matrices $A_j$ as well.
\end{rem}

We presume that there holds the following evident necessary condition 
 
\[ \sum {\rm Tr}(c_j)=0~,~{\rm resp.~}\prod \det (C_j)=1~.\] 
In terms of the eigenvalues $\lambda _{k,j}$  
(resp. $\sigma _{k,j}$) of the matrices from $c_j$ (resp. $C_j$) repeated 
with their multiplicities, this condition reads     

\[ \sum _{k=1}^n\sum _{j=1}^{p+1}\lambda _{k,j}=0~,~{\rm resp.~}
\prod _{k=1}^n\prod _{j=1}^{p+1}\sigma _{k,j}=1~.\]  
An equality of the kind
 
\[ \sum _{j=1}^{p+1}\sum _{k\in \Phi _j}\lambda _{k,j}=0~,~{\rm resp.~}
\prod _{j=1}^{p+1}\prod _{k\in \Phi _j}\sigma _{k,j}=1\]
is called a 
{\em non-genericity relation};  
the sets $\Phi _j$ contain one and the same number $<n$ of indices  
for all $j$. Eigenvalues satisfying none of these relations are called 
{\em generic}. Reducible  
$(p+1)$-tuples exist only for non-generic eigenvalues; indeed, the 
eigenvalues of each diagonal block of a block upper-triangular  
$(p+1)$-tuple satisfy some non-genericity relation.

\begin{defi}
Denote by $\{ J_j^n\}$  
a $(p+1)$-tuple of JNFs, $j=1$,$\ldots$, $p+1$. 
We say that the DSP is {\em solvable} (resp. that it is {\em weakly solvable} 
or, equivalently, that the weak DSP is solvable) for a given 
$\{ J_j^n\}$ and for given eigenvalues if there exists an 
irreducible $(p+1)$-tuple (resp. a $(p+1)$-tuple with a trivial centralizer) 
of matrices $M_j$ or of matrices $A_j$, with $J(M_j)=J_j^n$ or 
$J(A_j)=J_j^n$ and with the given 
eigenvalues. By definition, the DSP is solvable for $n=1$. Solvability of 
the DSP imlies its weak solvability, i.e. solvability of the weak DSP.
\end{defi}

For a given JNF $J^n=\{ b_{i,l}\}$ define its {\em corresponding} 
diagonal JNF ${J'}^n$. A diagonal JNF is  
a partition of $n$ defined by the multiplicities of the eigenvalues. 
For each $l$ $\{ b_{i,l}\}$ is a partition of $\sum _{i\in I_l}b_{i,l}$ and 
${J'}^n$ is the disjoint sum of the dual partitions. We say that two JNFs of 
one and the same size correspond to one another if they correspond to one and 
the same diagonal JNF.

\begin{prop}\label{rd}
1) One has $r(J^n)=r({J'}^n)$ and $d(J^n)=d({J'}^n)$.

2) To each diagonal JNF there corresponds a unique JNF with a single 
eigenvalue.
\end{prop}

\begin{ex}
To the JNF $\{ \{ 4,3,3\} ,\{ 3,2\} \}$ of size 15 
(two eigenvalues, with respectively 
three Jordan blocks, of sizes 4,3,3 and with two Jordan blocks, of sizes 3,2) 
there corresponds the diagonal JNF with multiplicities of the eigenvalues 
equal to 3,3,3,2,2,1,1. Indeed, the partition of $10$ dual to 4,3,3 
is 3,3,3,1; the partition of $5$ dual to 3,2 is 2,2,1. After this we 
arrange the multiplicities in decreasing order. 

To the two above JNFs there corresponds the JNF with a single eigenvalue with 
sizes of the Jordan blocks equal to 7,5,3. Indeed, 7,5,3 is the partition of 
15 dual to 3,3,3,2,2,1,1.
\end{ex}

For a given $\{ J_j^n\}$ with $n>1$, which satisfies condition 
$(\beta _n)$ and doesn't satisfy condition 

\[ (r_1+\ldots +r_{p+1})\geq 2n~~~~~~~~~~~~~~~~(\omega _n)\] 
set $n_1=r_1+\ldots +r_{p+1}-n$. Hence, $n_1<n$ and 
$n-n_1\leq n-r_j$. Define 
the $(p+1)$-tuple $\{ J_j^{n_1}\}$ as follows: to obtain the JNF 
$J_j^{n_1}$ 
from $J_j^n$ one chooses one of the eigenvalues of $J_j^n$ with 
greatest number $n-r_j$ of Jordan blocks, then decreases  
by 1 the sizes of the $n-n_1$ {\em smallest} Jordan blocks with this 
eigenvalue and deletes the Jordan blocks of size 0. 

\begin{defi}
The quantity $\kappa =2n^2-\sum _{j=1}^{p+1}d_j$ defined for a $(p+1)$-tuple 
of conjugacy classes is called the 
{\em index of rigidity}. It is introduced by N. Katz in \cite{Ka}. For 
irreducible representations it 
takes the values 2, 0, $-2$, $-4$, $\ldots$. Indeed, every conjugacy class is 
of even dimension and there holds condition $(\alpha _n)$. If for 
an irreducible $(p+1)$-tuple one has $\kappa =2$, then the $(p+1)$-tuple 
is called {\em rigid}. Such irreducible $(p+1)$-tuples are unique up to 
conjugacy, see \cite{Ka} and \cite{Si}. 
\end{defi}

\begin{lm}
The index of rigidity is invariant for the 
construction $\{ J_j^n\}\mapsto \{ J_j^{n_1}\}$.
\end{lm}

\begin{tm}\label{generic}
Let $n>1$. The DSP is solvable for the conjugacy classes $C_j$ or 
$c_j$ (with generic eigenvalues,  
defining the JNFs $J_j^n$ and satisfying condition  
$(\beta _n)$) if and only if either $\{ J_j^n\}$ satisfies 
condition $(\omega _n)$ or the construction 
$\{ J_j^n\}\mapsto \{ J_j^{n_1}\}$ iterated as long as it is defined 
stops at a $(p+1)$-tuple $\{ J_j^{n'}\}$ either with $n'=1$ or satisfying 
condition $(\omega _{n'})$.
\end{tm}

\begin{rems}
1) The conditions of the theorem are necessary for the weak solvability of the 
DSP for any eigenvalues. 

2) A posteriori one knows that the theorem does not 
depend on the choice(s) of eigenvalue(s) made when defining the 
construction $\{ J_j^n\}\mapsto \{ J_j^{n_1}\}$.
\end{rems}

\section{An example with index of rigidity equal to 2\protect\label{ex1}}

\subsection{Description of the example}

Denote by $J^*$, $J^{**}$ two quadruples of JNFs $J_j$ of size 4, 
$j=1,\ldots ,4$, in both of 
which $J_1$, $J_2$ and $J_3$ are diagonal, each with two eigenvalues of 
multiplicity 2; in $J^*$ the JNF $J_4$ is with a single eigenvalue to which 
there correspond three Jordan blocks, of sizes 2,1,1; in $J^{**}$ the JNF 
$J_4$ is diagonal, with two eigenvalues, of multiplicities 3 and 1. 
The JNFs $J_4$ from the two quadruples correspond to each other.

Hence, both $J^*$ and $J^{**}$ satisfy the conditions of 
Theorem~\ref{generic} (to be checked by the reader). They are both with 
index of rigidity 2. In both cases (of matrices $A_j$ or $M_j$) 
the quadruple $J^{**}$ 
admits generic eigenvalues and, hence, there exist irreducible quadruples of 
matrices $A_j$ or $M_j$ with such respective JNFs.

\begin{defi}
Suppose that the greatest common divisor of the multiplicities of all 
eigenvalues of the matrices $M_j$ or $A_j$ equals $q>1$. 
In the case of matrices $M_j$ denote by $\xi$ the product of 
all eigenvalues with multiplicities decreased $q$ times. Hence, $\xi$ is a 
root of unity of order $q$: $\xi =\exp (2\pi il/q)$, $l\in {\bf N}$. 
Denote by $m$ the greatest common divisor of $l$ and $q$. Hence, for $m>1$ 
the eigenvalues satisfy the 
non-genericity relation (called {\em basic}) their product with multiplicities 
divided $m$ times to equal 1. In the case of matrices $A_j$ the basic 
non-genericity relation is the sum of all eigenvalues with multiplicities 
decreased $q$ times to equal 0. Eigenvalues satisfying only the basic 
non-genericity relation and its corollaries are called 
{\em relatively generic}. 
\end{defi}

The quadruple $J^*$ does not admit generic but only relatively generic 
eigenvalues in the case of 
matrices $A_j$ because one has $q=2$. 

The quadruple $J^*$ admits generic eigenvalues in the case of matrices $M_j$. 
Indeed, such is the set of eigenvalues of the four matrices $(e,e^{-1})$, 
$(\sqrt{2},1/\sqrt{2})$, $(3,1/3)$, $i$. In this case $q=2$ and 
the product of all 
eigenvalues with multiplicities decreased twice equals $-1$. This is 
not a non-genericity relation. If the eigenvalue 
of the fourth matrix is changed from $i$ to $-1$, then the eigenvalues will 
not be generic -- their product when the multiplicities are decreased twice 
equals 1. This is the basic non-genericity relation. In this case 
the eigenvalues are relatively generic but not generic.

In our example we consider conjugacy classes $C_j$ defining the quadruple of 
JNFs $J^*$, with relatively generic but not generic eigenvalues. Observe that 
the expected dimension of ${\cal V}$ both in the case of $J^*$ and of $J^{**}$ 
equals $8+8+8+6-15=15$.

\subsection{The stratified structure of the variety ${\cal V}$ from the 
example}

The variety ${\cal V}$ from the example 
%
contains at least the following two 
strata denoted by ${\cal U}$ and ${\cal W}$. The stratum 
${\cal U}$ consists of all quadruples defining representations which are 
direct 
sums of two irreducible representations, i.e. up to conjugacy one has 
(for $(M_1,M_2,M_3,M_4)\in {\cal U}$)

\begin{equation}\label{direct} 
M_j=\left( \begin{array}{cc}N_j&0\\0&P_j\end{array}\right) ~,~N_j,P_j\in 
GL(2,{\bf C})
\end{equation}
where the matrices $N_j$ (resp. $P_j$) are diagonal for $j=1,2,3$. Their 
quadruples are with 
generic eigenvalues and for 
$j=4$ the eigenvalues equal $-1$, $P_4$ is conjugate to a Jordan block of 
size 2 while $N_4$ is scalar. The existence of irreducible quadruples of 
matrices $N_j$ and $P_j$ is guaranteed by Theorem~\ref{generic}.

\begin{rem}
The matrices $N_j$ (resp. $P_j$) define an irreducible rigid representation 
(resp. an irreducible representation of zero index of rigidity).
\end{rem}

\begin{prop}\label{N_jP_j}
1) The variety of matrices $N_j$ (resp. $P_j$) as above is smooth, irreducible 
and of dimension 3 (resp. 5).

2) The variety of quadruples of diagonalizable matrices 
$M_j\in GL(2,{\bf C})$ each with two 
distinct eigenvalues (the eigenvalues of the quadruple being generic) is 
smooth, irreducible and of dimension 5.
\end{prop}
 
All propositions from this subsection are proved in Section~\ref{proofofUW}.

The stratum ${\cal W}$ consists of all quadruples defining semi-direct sums of 
two equivalent rigid representations. Up to conjugacy one has 
(for $(M_1,M_2,M_3,M_4)\in {\cal W}$)

\begin{equation}\label{semidirect} 
M_j=\left( \begin{array}{cc}N_j&R_j\\0&N_j\end{array}\right) ~,~N_j,R_j\in 
GL(2,{\bf C})
\end{equation}
with $N_j$ as above. The blocks $R_j$ are such that for $j=1,2,3$ the matrices 
$M_j$ are diagonalizable while $M_4$ has JNF $J_4$ (i.e. rk$R_4=1$). 

The absence of other possible types of representations is guaranteed by the 
following theorem which follows from Theorem 1.1.2 from \cite{Ka}. The theorem 
and its proof were suggested by Ofer Gabber.

\begin{tm}\label{Gabber} 
For fixed conjugacy classes with index of rigidity 2 there cannot coexist 
irreducible and reducible $(p+1)$-tuples of matrices $M_j$.
\end{tm}

The theorem is proved in the Appendix. It follows from the theorem that there 
can exist only reducible quadruples of matrices $M_j$ in the example under 
consideration. 

\begin{prop}\label{only}
One has ${\cal V}={\cal U}\cup {\cal W}$.
\end{prop}

\begin{prop}\label{R_4}
1) In a quadruple (\ref{semidirect}) the matrix $R_4$ is nilpotent of rank 1 
and for $j=1,2,3$ one has $R_j=[N_j,Z_j]$ with $Z_j\in sl(2,{\bf C})$. 

2) If the matrices $N_1$, $N_2$, $N_3$ are fixed, then for every nilpotent 
rank 1 matrix $R_4$ there exists a quadruple of matrices (\ref{semidirect}). 
\end{prop}

\begin{prop}\label{centralizer}
The centralizers in $SL(4,{\bf C})$ of the quadruples (\ref{direct}) and 
(\ref{semidirect}) are both of dimension 1. They consist respectively 
of the matrices
\[ \left( \begin{array}{cc}\alpha I&0\\0&\pm \alpha ^{-1}I\end{array}\right) ~
{\rm and~}\left( \begin{array}{cc}\delta I&\beta I\\0&\delta I
\end{array}\right) ~,~\alpha \in {\bf C}^*~,~\beta \in {\bf C}~,~\delta ^4=1\]
\end{prop}

\begin{prop}\label{UW}
The stratum ${\cal W}$ belongs to the closure of the stratum ${\cal U}$. 
\end{prop}

\begin{prop}\label{dimW}
The stratum ${\cal W}$ is an irreducible smooth variety of dimension 15. 
\end{prop}

\begin{prop}\label{dimU}
The stratum ${\cal U}$ is an irreducible smooth variety of dimension 16.
\end{prop}

\begin{rem}
The closure of the variety ${\cal W}$ (hence, the one of ${\cal U}$ as well) 
contains the variety ${\cal Y}$ 
of quadruples which up to conjugacy are 
of the form (\ref{semidirect}) with $R_j=0$ for all $j$. For such quadruples 

1) the matrix $M_4$ is scalar; 

2) they define direct sums of two equivalent irreducible 
rigid representations. 

There exist no irreducible such quadruples of matrices 
$M_j$ or $A_j$ because the conditions of Theorem~\ref{generic} are not 
fulfilled (neither the necessary condition $(\alpha _n)$). 
\end{rem}

\begin{prop}\label{dimY}
The variety ${\cal Y}$ is smooth and irreducible. One has dim${\cal Y}=12$. 
\end{prop}

\section{Proofs of the propositions\protect\label{proofofUW}}

{\bf Proof of Proposition~\ref{N_jP_j}:}

$1^0$. The variety of quadruples of matrices $N_j$ is obtained by 
conjugating one such quadruple by matrices from $SL(2,{\bf C})$ (indeed, 
rigid $(p+1)$-tuples are unique up to conjugacy, see \cite{Ka} and \cite{Si}). 
This proves the connectedness. The smoothness and the dimension follow 
from Proposition~\ref{smooth}. 

$2^0$. Denote by $C_j^*$ the conjugacy class of the matrix $P_j$. 
Prove that the variety $\Pi$ of quadruples of matrices $P_j$ is 
connected. 
Denote by $\delta$ the product det$P_1$det$P_2$. By varying the matrices $P_1$ 
and $P_2$ (resp. $P_3$ and $P_4$) one can obtain as their product $P_1P_2$ 
(resp. as $P_4^{-1}P_3^{-1}$) any matrix from the 
set $\Delta (\delta )$ of $2\times 2$-matrices with determinant 
equal to $\delta$. The set $\Delta (\delta )$ being connected so is the 
variety $\Pi$ because 
$\Pi =\{ (P_1,P_2,P_3,P_4)|P_j\in C_j^*,P_1P_2=P_4^{-1}P_3^{-1}\}$.   

$3^0$. The eigenvalues of the matrices $P_j$ being generic the variety 
$\Pi$ contains no reducible quadruples. Hence, the variety $\Pi$ is smooth, 
one has dim$\Pi =5$, see Proposition~\ref{smooth}. 

$4^0$. Part 2) is proved by analogy with $2^0$ and $3^0$.\hspace{1cm}$\Box$\\ 

{\bf Proof of Proposition~\ref{only}:} 

$1^0$. A quadruple from ${\cal V}$ is block upper-triangular up to conjugacy. 
The eigenvalues being relatively generic the diagonal blocks can be only of 
size 2 and the restrictions of the matrices $M_j$ to them can be with 
conjugacy classes like in the cases of quadruples of matrices $N_j$ or $P_j$. 

$2^0$. Show that if one of the diagonal blocks is a quadruple 
of matrices $N_j$ and the other one of matrices $P_j$, then this is a direct 
sum conjugate to a quadruple (\ref{direct}). Indeed, 
for the representations $P$ and $N$ 
defined by the quadruples of matrices $P_j$ and $N_j$ one has 
Ext$^1(P,N)=$Ext$^1(N,P)=0$ (to be checked directly). This implies that a 
block upper-triangular quadruple of matrices $M_j$ with diagonal 
blocks $N_j$ and $P_j$ is conjugate to its restriction to the two 
diagonal blocks, i.e. the quadruple is a point 
from ${\cal U}$. On the other hand, if both diagonal blocks equal $N_j$, then 
the quadruple is like in (\ref{semidirect}).  

Hence, only quadruples like the ones from ${\cal U}$ and ${\cal W}$ can 
exist in ${\cal V}$.\hspace{1cm}$\Box$\\   

{\bf Proof of Proposition~\ref{R_4}:}

$1^0$. The blocks $R_1$, $R_2$, $R_3$ must be of the form $R_j=[N_j,Z_j]$ for 
some matrices $Z_j\in gl(n,{\bf C})$. Indeed, it suffices to prove this in the 
assumption that $N_j$ is diagonal: 
$N_j=\left( \begin{array}{cc}\lambda &0\\0&\mu \end{array}\right)$, 
$\lambda \neq \mu$. Set  
$R_j=\left( \begin{array}{cc}g&h\\f&s\end{array}\right)$. One must have 
$g=s=0$, otherwise $M_j$ will not be diagonalizable. But then 
$R_j=[N_j,Z_j]$ with 
$Z_j=\left( \begin{array}{cc}0&h/(\lambda -\mu )\\f/(\mu -\lambda )&0
\end{array}\right)$. 

On the other hand, if for $j=1,2,3$ one has $R_j=[N_j,Z_j]$, then the 
matrices $M_1$, $M_2$, $M_3$ have the necessary 
JNFs -- one has 
$M_j=\left( \begin{array}{cc}I&Z_j\\0&I\end{array}\right) ^{-1} 
\left( \begin{array}{cc}N_j&0\\0&N_j\end{array}\right) 
\left( \begin{array}{cc}I&Z_j\\0&I\end{array}\right)$ .

$2^0$. If one has rk$R_4=0$, then $R_4=0$ and $M_4$ must be scalar, i.e. 
$M_4\not\in C_4$. If rk$R_4=2$, then rk$(M_4+I)=2$ and again 
$M_4\not\in C_4$. Hence, rk$R_4=1$. This leaves two possibilities -- either 
$R_4$ has two distinct eigenvalues one of which is 0 or it is nilpotent.

$3^0$. The condition $M_1\ldots M_4=I$ restricted to the right upper block 
and to each of the diagonal blocks reads respectively

\[ R_1N_2N_3N_4+N_1R_2N_3N_4+N_1N_2R_3N_4+N_1N_2N_3R_4=0~,~N_1N_2N_3=-I~.\]
Hence, the first of these two equalities takes the form 

\[ -R_1-N_1R_2(N_1)^{-1}-(N_1N_2)R_3(N_1N_2)^{-1}-R_4=0~.\]
As $R_j=[N_j,Z_j]$, $j=1,2,3$, see $1^0$, one has 

\[ {\rm tr}R_1={\rm tr}R_2={\rm tr}(N_1R_2(N_1)^{-1})={\rm tr}R_3=
{\rm tr}((N_1N_2)R_3(N_1N_2)^{-1})=0~.\]
Hence, tr$R_4=0$. This means that $R_4$ is nilpotent, 
of rank 1. This proves 1).

$4^0$. To prove 2) one has to recall that $R_j=[N_j,Z_j]$ for $j=1,2,3$, 
see $1^0$,  
and that each matrix from $sl(2,{\bf C})$ can be represented as 
$\sum _{j=1}^3[N_j,Z_j]$, see Proposition~\ref{commut}. Hence, for every 
nilpotent $R_4$ one can find matrices $Z_j$ such that for $j=1,2,3$ one has 
$R_j=[N_j,Z_j]$, i.e. $M_j\in C_j$ and $M_1M_2M_3M_4=I$.\hspace{1cm}$\Box$\\

{\bf Proof of Proposition~\ref{centralizer}:}

$1^0$. Denote by $F=\left( \begin{array}{cc}U&V\\W&Y\end{array}\right)$ a 
matrix from the centralizer of the quadruple. In the case of a quadruple 
(\ref{direct}) the commutation relations read:

\[ [U,N_j]=[Y,P_j]=0~,~N_jV=VP_j~,~WN_j=P_jW~.\]
The representations defined by the matrices $N_j$ and $P_j$ being 
non-equivalent, these relations imply $V=W=0$. The irreducibility of the 
quadruples of matrices $N_j$ and $P_j$ and Schur's lemma imply that $U$ and 
$Y$ are scalar. Hence, $U=\alpha I$, $Y=\xi I$ with $\alpha ^2\xi ^2=1$, 
i.e. $\xi =\pm \alpha ^{-1}$. 

$2^0$. In the case of a quadruple (\ref{semidirect}) the matrix algebra 
${\cal A}$ generated by the matrices $M_j$ contains the matrix $M_4+I$ and 
its left 
and right products by matrices from the algebra ${\cal B}$ 
generated by $M_1$, $M_2$ and $M_3$. As ${\cal B}$ contains matrices of the 
form $\left( \begin{array}{cc}T&\ast \\0&T\end{array}\right)$ for any 
$T\in gl(2,{\bf C})$ (the Burnside theorem), the algebra ${\cal A}$ contains 
all matrices of the form $\left( \begin{array}{cc}0&Q\\0&0\end{array}\right)$ 
with $Q\in gl(2,{\bf C})$. 

The commutation relations 
imply that $WQ=0$, hence, $W=0$, and $UQ=QY$ for any $Q$, i.e. $U=Y=\delta I$. 
Finally, one has $[N_j,V]=0$ which implies that 
$V=\beta I$ (use Schur's lemma). 

One must have $\delta ^4=1$ because $F\in SL(4,{\bf C})$.\hspace{1cm}$\Box$\\

{\bf Proof of Proposition~\ref{UW}:}

$1^0$. One can deform the matrices $M_j$ from a quadruple from ${\cal W}$ as 
follows. The deformation parameter is denoted by 
$\varepsilon \in ({\bf C},0)$ and the deformed matrices by $M_j'$. 
Assume that $N_4=-I$, $R_4=\left( \begin{array}{cc}0&1\\0&0\end{array}\right)$ 
(one can achieve this by conjugation of the quadruple 
with a block-diagonal matrix). Set 
$M_4'=M_4+\varepsilon (E_{1,2}+w(\varepsilon )E_{1,3})$; the matrix $E_{k,j}$ 
by definition 
has a single non-zero entry equal to 1 in position $(k,j)$; $w(\varepsilon )$ 
is an unknown germ of an analytic function. 

$2^0$. For $j=1,2,3$ set 
$M_j'=(I+\varepsilon X_j(\varepsilon ))^{-1}M_j
(I+\varepsilon X_j(\varepsilon ))$ 
where $X_j=\left( \begin{array}{cc}U_j&V_j\\0&0\end{array}\right)$. 
Set $X_j(0)=X_j^0$.
One must have $M_1'M_2'M_3'M_4'=I$ which in 
first approximation w.r.t. $\varepsilon$ reads 
\begin{equation}\label{FA} 
[M_1,X_1^0]M_2M_3M_4+M_1[M_2,X_2^0]M_3M_4+M_1M_2[M_3,X_3^0]M_4+
M_1M_2M_3(E_{1,2}+w(0)E_{1,3})=0
\end{equation}

$3^0$. Set $U_j^0=U_j(0)$, $U^0=(U_1^0,U_2^0,U_3^0)$, $V_j^0=V_j(0)$, 
$V^0=(V_1^0,V_2^0,V_3^0)$, $w^0=w(0)$. Equation (\ref{FA}) 
restricted to the left upper block reads:

\[ {\cal G}(U^0):=[N_1,U_1^0]N_2N_3+N_1[N_2,U_2^0]N_3+
N_1N_2[N_3,U_3^0]=N_1N_2N_3E_{1,2}\]
(because $N_4=-I$). Making use of $N_1N_2N_3=-I$ one finds 

\begin{equation}\label{N} 
[N_1,U_1^0N_1^{-1}]+[N_1N_2N_1^{-1},N_1U_2^0N_2^{-1}N_1^{-1}]+
[N_3,N_3^{-1}U_3^0]=E_{1,2}
\end{equation}  
The triple of matrices $N_1$, $N_2$, $N_3$ is irreducible, hence, so is the 
triple $N_1$, $N_1N_2N_1^{-1}$, $N_3$. By Proposition~\ref{commut}, one can 
find matrices $U_j^0$ satisfying equation (\ref{N}).

$4^0$. Equation (\ref{FA}) restricted to the right lower block is of the 
form $0=0$, i.e. it gives no condition at all upon $U_j^0$, $V_j^0$ and $w^0$. 
Its restriction to the right upper block reads:

\begin{equation}\label{V}
{\cal F}(V^0,U^0,w^0):=
{\cal G}(V^0)+{\cal H}(U^0)-w^0E_{1,3}=0
\end{equation}
where ${\cal H}$ is some linear form in the entries of the matrices $U_j^0$. 
Hence, if $U_j^0$ are found such that (\ref{N}) holds, then one can find 
$w^0$ such that tr$({\cal H}(U_1^0,U_2^0,U_3^0))=w^0$. After this one can 
find matrices $V_j^0$ such that (\ref{V}) hold.

$5^0$. The map $(U^0,V^0,w^0)\mapsto ({\cal G}(U^0),{\cal F}(V^0,U^0,w^0))$ is 
surjective onto the space of $2\times 4$-matrices. By the implicit function 
theorem one can find germs of matrices $U_j$, $V_j$ and a germ of a 
function $w$ holomorphic in $\varepsilon$ at 0 such that $M_1'\ldots M_4'=I$. 

Fix $\varepsilon \neq 0$. 
The quadruple of matrices $M_j'$ is block upper-triangular with diagonal 
blocks having the properties of $P_j$ and $N_j$ ($P_j$ is above). Moreover, 
each of the matrices $M_j'$ is conjugate to the block-diagonal matrix whose 
restriction to the two diagonal blocks is the same as the one of $M_j'$ 
(to be checked directly). By Proposition~\ref{only}, up to conjugacy the 
quadruple of matrices is 
like the one from (\ref{direct}).\hspace{1cm}$\Box$\\

{\bf Proof of Proposition \ref{dimW}:}

$1^0$. Prove the irreducibility.  
The variety ${\cal W}$ is obtained by 
conjugating with matrices from $SL(4,{\bf C})$ the quadruples of matrices of 
the form (\ref{semidirect}) with $R_4$ nilpotent of rank 1. The orbit of 
$R_4$ is an irreducible variety which implies the irreducibility of 
${\cal W}$.

$2^0$. Fix the blocks $N_j$ of a 
quadruple (\ref{semidirect}). The variety ${\cal S}$ 
of such quadruples defined modulo 
conjugacy is of dimension 1. Indeed, the orbit of $R_4$ is of dimension 2. 
The only conjugations that preserve the form of the quadruple and its 
restrictions to the two diagonal blocks are with matrices of the form 
$\left( \begin{array}{cc}aI&V\\0&bI\end{array}\right)$, $ab\neq 0$, 
$V\in gl(2,{\bf C})$; this 
is proved in $4^0$. If one 
requires the matrix to be from $SL(4,{\bf C})$, this means that $b=\pm 1/a$ 
and factoring out these conjugations decreases the dimension by 1. Indeed, 
such a conjugation changes $R_4$ to $bR_4/a$, 
the presence of $V$ does not affect the block $R_4$. 

$3^0$. To obtain the variety ${\cal H}$ of all quadruples defining semi-direct 
sums like 
(\ref{semidirect}) one has to conjugate the quadruples from ${\cal S}$ by 
matrices from $SL(4,{\bf C})$. This increases the dimension by 14 (not by 15 
because the centralizer of such a quadruple is non-trivial, of dimension 1, 
see Proposition~\ref{centralizer}). Hence, dim${\cal H}=15$.

$4^0$. Denote by $G$ a matrix the conjugation with which preserves 
the block upper-triangular form of the quadruple and the blocks $N_j$. 
If $G=\left( \begin{array}{cc}U&V\\W&Y\end{array}\right)$, 
then the condition the quadruple to remain block upper-triangular implies that 
$[W,N_j]=0$, i.e. $W=hI$. The condition the diagonal blocks of $M_4$ to 
remain the same implies $[N_4,Y]-WR_4=R_4W+[N_4,U]=0$. As $N_4=-I$, one has 
$[N_4,Y]=[N_4,U]=0$, i.e. $W=0$. 

The conditions $[N_j,U]=[N_j,Y]=0$ imply that 
$U=aI$, $Y=bI$.\hspace{1cm}$\Box$\\

{\bf Proof of Proposition \ref{dimU}:}

$1^0$. The varieties of quadruples of matrices $N_j$ or $P_j$, see  
Proposition~\ref{N_jP_j}, are smooth, irreducible and of dimensions 
respectively 
3 and 5. Hence, the 
variety ${\cal P}$ of quadruples of matrices $M_j$ like in 
(\ref{direct}) is smooth, irreducible and of dimension 8. 

$2^0$. The variety ${\cal U}$ is of dimension $8+15-7=16$. Here 
``8'' stands for ``dim${\cal P}$'', ``15'' stands for 
``dim$SL(4,{\bf C})$'' and 7 is the dimension of the subgroup of 
$SL(4,{\bf C})$ of block-diagonal matrices with blocks $2\times 2$ 
conjugation with which preserves the block-diagonal 
form of quadruple (\ref{direct}) (infinitesimal conjugations {\em only} with 
such matrices preserve the block-diagonal form of quadruple (\ref{direct})); 
this subgroup contains the centralizer of quadruple (\ref{direct}), 
see Proposition~\ref{centralizer}. \hspace{1cm}$\Box$\\ 

{\bf Proof of Proposition \ref{dimY}:} 

The variety ${\cal Y}$ is the 
orbit of one quadruple of the form (\ref{semidirect}) with $R_j=0$, 
$j=1,\ldots ,4$, under conjugation by $SL(4,{\bf C})$ 
(recall that the matrices $N_j$ define a rigid representation, i.e. unique up 
to conjugacy). Hence, ${\cal Y}$ is irreducible and smooth.

To obtain dim${\cal Y}$ one has to subtract from 
$15=$dim$SL(4,{\bf C})$ the dimension of the centralizer in $SL(4,{\bf C})$ 
of the above quadruple. The latter equals 3 -- the centralizer 
is the set of all matrices of the form 
$\left( \begin{array}{cc}\alpha I&\beta I\\ \delta I&\eta I
\end{array}\right)$ with $\alpha \eta -\delta \beta=\pm 1$.\hspace{1cm}$\Box$

\section{Another example with index of rigidity 2\protect\label{ex2}}

Consider the variety ${\cal V}$ in the case when $p=2$, $n=4$, the three 
conjugacy classes are diagonalizable and have eigenvalues $(a,a,b,c)$, 
$(f,f,g,h)$ and $(u,u,v,w)$ (different letters denote different eigenvalues). 
The index of rigidity equals 2 (to be checked directly). 

The eigenvalues are presumed to 
satisfy the only non-genericity relation $abfguv=1$. Hence, 
for such conjugacy classes there exist irreducible triples of 
diagonalizable matrices 
$L_j\in gl(2,{\bf C})$ (resp. $B_j\in gl(2,{\bf C})$) with eigenvalues 
$(a,b)$; $(f,g)$; $(u,v)$ (resp. $(a,c)$; 
$(f,h)$; $(u,w)$) such that $L_1L_2L_3=I$ (resp. $B_1B_2B_3=I$). This follows 
from Theorem~\ref{generic}. Hence, there 
exist triples of block-diagonal matrices $M_j$ with diagonal blocks equal 
to $L_j$ and $B_j$. Denote by ${\cal D}$ the variety of such triples.
By Theorem~\ref{Gabber}, irreducible triples of matrices $M_j$ do not exist. 

There do exist, however, triples with trivial centralizers which are 
block upper-triangular: $M_j=\left( \begin{array}{cc}L_j&T_j\\0&B_j
\end{array}\right)$ where $T_j=L_jY_j-Y_jB_j$ for some 
$Y_j\in gl(2,{\bf C})$ because $M_j$ is conjugate to 
$\left( \begin{array}{cc}L_j&0\\0&B_j
\end{array}\right)$. The condition $M_1M_2M_3=I$ restricted to the right upper 
block reads: 

\[ T_1B_2B_3+L_1T_2B_3+L_1L_2T_3=0~~~~~~(*)\]  
Thus the triple of matrices $T_j$ belongs to the space 

\[ {\cal T}=\{ (T_1,T_2,T_3)~|~T_j=L_jY_j-Y_jB_j~,~Y_j\in gl(2,{\bf C})~,~
T_1B_2B_3+L_1T_2B_3+L_1L_2T_3=0\} ~.\]

{\em One has dim${\cal T}=5$.} 

Indeed, the conditions $T_j=L_jY_j-Y_jB_j$ imply that each matrix $T_j$ 
belongs to the image of the map $(.)\mapsto L_j(.)-(.)B_j$ which is a 
subspace of $gl(2,{\bf C})$ of dimension 3. Condition (*) 
is equivalent to four linearly independent equalities (we let the reader 
prove their linear independence using the non-equivalence of the 
representations defined by the matrices $L_j$ and $B_j$). 

Consider the space 

\[ {\cal Q}=\{ (T_1,T_2,T_3)~|~T_j=L_jY-YB_j~,~Y\in gl(2,{\bf C})\} ~.\]
For such matrices $T_j$ there holds $(*)$, therefore 
${\cal Q}\subset {\cal T}$. 
The space ${\cal Q}$ is the space of right upper blocks of triples of block 
upper-triangular matrices $M_j$ which are obtained from block-diagonal ones 
from ${\cal D}$ by conjugation with matrices of the form 
$\left( \begin{array}{cc}I&Y\\0&I\end{array}\right)$.  

{\em One has dim${\cal Q}=4$.}

Indeed, for no matrix from $gl(2,{\bf C})$ does one have $L_jY-YB_j=0$ for 
$j=1,2,3$ because the triples of matrices $L_j$ and $B_j$ define 
non-equivalent representations.

Hence, dim$({\cal T}/{\cal Q})=1$. Choose the triple of 
matrices $Y_j$ to span the factorspace $({\cal T}/{\cal Q})$. Hence, the 
centralizer ${\cal Z}$ of the 
triple of matrices $M_j$ will be trivial. Indeed, let 
$Z=\left( \begin{array}{cc}P&Q\\R&S\end{array}\right) \in {\cal Z}$. 
Hence, $RL_j=B_jR$ for $j=1,2,3$ (commutation relations restricted to 
the left lower block), i.e. $R=0$ because the matrices $L_j$ and $B_j$ 
define non-equivalent representations. 

One must have $[P,L_j]=[S,B_j]=0$ (commutation relations restricted to the 
diagonal blocks), 
i.e. $P=aI$, $B=bI$. But then one must have (commutation relations 
restricted to the right upper block) $(a-b)T_j=L_jQ-QB_j$ which 
means that $a=b$ (otherwise $(T_1,T_2,T_3)\in {\cal Q}$), hence, 
$L_jQ-QB_j=0$ for $j=1,2,3$, i.e. $Q=0$. Hence, $Z=aI$. 

\begin{rems}
1) It is clear that the variety ${\cal D}$ belongs to the closure of 
${\cal V}\backslash {\cal D}$ -- the triple of matrices 
$M_j=\left( \begin{array}{cc}L_j&\varepsilon T_j\\0&B_j\end{array}\right)$ 
belongs to ${\cal V}\backslash {\cal D}$ for $\varepsilon \neq 0$, for 
$\varepsilon =0$ it belongs to ${\cal D}$.

2) The variety ${\cal V}$ is connected, hence, irreducible. This follows from 
$({\cal T}/{\cal Q})$ being a linear space (${\cal V}$ is obtained by 
conjugating block upper-triangular triples with 
$(T_1,T_2,T_3)\in ({\cal T}/{\cal Q})$ and with fixed diagonal blocks by 
matrices from $SL(4,{\bf C})$).
\end{rems}

\section{A third example with index of rigidity 2\protect\label{ex3}}

Let $n=4$, $p=2$. Use the notation from the previous section. Define the 
conjugacy classes $C_j$ as follows: their eigenvalues equal $(a,a,b,b)$, 
$(f,f,g,g)$, $(u,u,v,v)$, the eigenvalues are relatively generic but not 
generic (one has $abfguv=1$). To each of the eigenvalues $a$, 
$b$ and $f$ there corresponds 
a single Jordan block of size 2, to each of the eigenvalues $g$, $u$, $v$ 
there correspond two Jordan blocks of size 1. Hence, the index of rigidity 
equals 2.

The variety ${\cal V}$ contains triples of matrices which up to conjugacy are 
block upper-triangular with two diagonal blocks equal to $L_j$, see their 
definition in the  
previous section. By Theorem~\ref{Gabber}, ${\cal V}$ contains no irreducible 
triples. Hence, it contains none with trivial centralizer either because 
the matrices $M_j$ from any such block upper-triangular 
triple commute with the matrix $E_{1,3}+E_{2,4}$; on the other hand, if a 
triple of matrices $M_j\in C_j$ is conjugated to a block upper-triangular 
form, then the diagonal blocks are of size 2 and up to conjugacy they 
equal $L_j$ -- this follows from the choice of the eigenvalues.  

\begin{prop}\label{dimV}
One has dim${\cal V}=15$ which is the expected dimension.
\end{prop}

\begin{rems}
The closure of the variety ${\cal V}$ contains the varieties in which at 
least one of the two Jordan normal forms $J(M_1)$ and $J(M_2)$ contains 
instead of some Jordan block(s) of size 2 two Jordan blocks of size 1. We 
leave the details for the reader. One can prove that ${\cal V}$ is 
irreducible.  
\end{rems}

{\bf Proof of Proposition~\ref{dimV}:}

$1^0$. Suppose that one has $M_j=\left( \begin{array}{cc}L_j&T_j\\0&L_j
\end{array}\right)$ with $L_1=$diag$(a,b)$, $T_1=$diag$(1,1)$. Fix $L_2$ 
and $L_3$. Then the couple of blocks $(T_2,T_3)$ belongs to a space  
of dimension 1. 

Indeed, one has $T_3=[L_3,Z_3]$ in order $M_3$ to be 
diagonalizable and the dimension of the image of the map 
$Z_3\mapsto [L_3,Z_3]$ in $gl(2,{\bf C})$ equals 2. 

The block $T_2$ belongs to 
an affine space of dimension 2. Indeed, one has $T_2=S+[L_2,Z_2]$,  
where the dimension of the image of the map $Z_2\mapsto [L_2,Z_2]$ equals 2 
and the matrix $S$ is defined as follows. Set $L_2=H^{-1}$diag$(f,g)H$. 
Then $S=\xi H^{-1}E_{1,3}H$ where $\xi$ satisfies the condition 

\[ {\rm tr}(L_2L_3+L_1SL_3)=0~~~~~~(**)\]
(If by chance this condition gives $\xi =0$, then one 
has to choose two diagonal entries of $T_1$ other than $(1,1)$ so that 
$\xi \neq 0$, otherwise $M_2$ will be diagonalizable.)

$2^0$. The coefficient $\xi$ satisfies condition (**) for the following 
reason. 
The condition $M_1M_2M_3=I$ implies that 
${\cal H}:=(T_1L_2L_3+L_1T_2L_3+L_1L_2T_3)=0$. In particular, 
tr${\cal H}=0$. As 

\[ L_1L_2L_3=I~,~T_1=I~,~ 
{\rm tr}(L_1L_2T_3)={\rm tr}(L_1L_2L_3Z_3-L_1L_2Z_3L_3)={\rm tr}
(Z_3-L_3^{-1}Z_3L_3)=0\] 
and tr$(L_1[L_2,Z_2]L_3)=$tr$(L_3^{-1}Z_2L_3-L_1Z_2L_1^{-1})=0$, one has 
tr$(L_2L_3+L_1SL_3)=0$. 

$3^0$. From the dimension $2+2$ of the space to which the couple $(T_2,T_3)$ 
belongs one has to subtract 3 because the equation ${\cal H}=0$ (after 
one has chosen $\xi$ so that tr${\cal H}=0$) imposes 3 conditions.

$4^0$. The centralizer ${\cal Z}$ 
of the triple of matrices $M_j$ in $SL(4,{\bf C})$ 
is generated by the 
matrix $E_{1,3}+E_{2,4}$. Moreover, any matrix from $SL(4,{\bf C})$ the 
conjugation with which preserves the form of the triple belongs 
to ${\cal Z}$. This can be proved by a direct computation 
which we leave for the reader. 

$5^0$. To find the dimension of ${\cal V}$ one has to conjugate the 
block upper-triangular triples 
from $1^0$ whose variety is of dimension 1 by matrices from 
$SL(n,{\bf C})/{\cal Z}$. The latter variety is of dimension 14. Hence, 
dim${\cal V}=15$.\hspace{1cm}$\Box$

\section{A fourth example with index of rigidity 2\protect\label{ex4}}

Let $n=p=3$ and let the conjugacy classes $C_j$ define diagonal but 
non-scalar JNFs the eigenvalues being equal respectively to $(a,1,1)$, 
$(b,1,1)$, $(c,1,1)$, $(d,1,1)$, with $abcd=1$. Hence, the index of 
rigidity is 0. There exist reducible such 
quadruples of matrices $M_j$ with trivial centralizers. Example:

\[ M_1=\left( \begin{array}{ccc}a&0&0\\0&1&0\\0&0&1
\end{array}\right) ~,~M_2=\left( \begin{array}{ccc}b&1&0\\0&1&0\\0&0&1
\end{array}\right) \]

\[ M_3=\left( \begin{array}{ccc}c&0&1\\0&1&0\\0&0&1
\end{array}\right) ~,~M_4=\left( \begin{array}{ccc}d&-1/bc&-1/c\\0&1&0\\0&0&1
\end{array}\right) \]
(the reader is invited to check the triviality of the centralizer 
oneself). Denote by ${\cal T}$ 
the stratum of ${\cal V}$ of quadruples with trivial centralizers. Hence, 
dim${\cal T}=8$ (Proposition~\ref{smooth}). By Theorem~\ref{Gabber}, there 
exist no irreducible quadruples of matrices $M_j\in C_j$.

On the other hand, there exist quadruples defining direct sums of an 
irreducible representation of rank 2 and of a one-dimensional one. Example:
  
\[ M_1=\left( \begin{array}{ccc}a&1&0\\0&1&0\\0&0&1
\end{array}\right) ~,~M_2=\left( \begin{array}{ccc}b&-1/a&0\\0&1&0\\0&0&1
\end{array}\right) \]

\[ M_3=\left( \begin{array}{ccc}c&0&0\\-1/d&1&0\\0&0&1
\end{array}\right) ~,~M_4=\left( \begin{array}{ccc}d&0&0\\1&1&0\\0&0&1
\end{array}\right) \]
Denote by ${\cal S}$ the stratum of ${\cal V}$ of quadruples defining such 
direct sums. 

{\em One has dim${\cal S}=9$.}  

Indeed, the subvariety 
${\cal S}'\subset {\cal S}$ 
of block-diagonal such quadruples is of dimension 5 
(Proposition~\ref{smooth}). Hence, ${\cal S}$ is obtained from ${\cal S}'$ by 
conjugating with matrices from $SL(3,{\bf C})$ (dim$SL(3,{\bf C})=8$) 
and one has to factor out the 
conjugation with block-diagonal matrices whose subgroup is of dimension 4. 
Thus dim${\cal S}=5+8-4=9$. 

\begin{rems}
1) Both strata ${\cal S}$ and ${\cal T}$ contain in their closures the variety 
of quadruples which are diagonal up to conjugacy, also the ones of quadruples 
defining direct sums of the one-dimensional representation 1, 1, 1, 1 with 
the semi-direct sums of the 
representations 1, 1, 1, 1 and $a$, $b$, $c$, $d$. 

2) The stratum ${\cal T}$ {\em does not} lie in the closure of the stratum 
${\cal S}$ (triviality of the centralizer is an ``open'' property). 

3) One can show that at every point of ${\cal V}$ one has dim${\cal V}\leq 9$. 
\end{rems}

\section{An example with zero index of rigidity\protect\label{ex5}}

By Theorem~\ref{generic}, there exist irreducible quadruples of matrices 
$A_j$ or $M_j$ of size 2 in which each matrix has two distinct eigenvalues and 
the eigenvalues are generic. 
For such quadruples the index of rigidity equals 0 (to be checked directly). 

Consider a quadruple of matrices (say, $M_j$; for matrices $A_j$ one can give 
a similar example) of the form 
  
\[ M_j=\left( \begin{array}{cc}B_j&0\\0&G_j\end{array}\right)\]
where each of the quadruples of matrices $B_j$ and $G_j$ is like above, with 
generic eigenvalues. Moreover, for each $j$ the eigenvalues of $B_j$ and $C_j$ 
are the same but the quadruples of matrices $B_j$ and $G_j$ define 
non-equivalent representations. To choose them such is possible because the 
quadruples are not rigid. 

Compute the dimension of the variety ${\cal M}$ 
of such quadruples of matrices $M_j$. The varieties ${\cal B}$ and ${\cal G}$ 
of quadruples of $2\times 2$-matrices $B_j$ or $G_j$ are both of 
dimension 5 (see part 2) of Proposition~\ref{N_jP_j}). 

Hence, dim${\cal M}=10$. The variety ${\cal N}$ 
of quadruples of matrices $M_j$ 
defining a direct sum of two representations of rank 2 with the properties 
of ${\cal B}$ and ${\cal G}$ is obtained by conjugating the 
quadruples from ${\cal M}$ by matrices from $SL(4,{\bf C})$. Infinitesimal 
conjugation 
by block-diagonal 
matrices from $SL(4,{\bf C})$ with two diagonal blocks of size 2 and only by 
such matrices preserves ${\cal M}$ (their 
subgroup is of dimension 7 in $SL(4,{\bf C})$). Hence, 
dim${\cal N}=10+15-7=18$ where $15=$dim$SL(4,{\bf C})$.   

The expected dimension of the variety ${\cal N}$ equals 17, see 
Proposition~\ref{smooth}.   
In a subsequent paper the author intends to prove that for zero 
index of rigidity and for relatively generic but not generic eigenvalues the 
Deligne-Simpson problem is not weakly solvable. Hence, in the above example 
one has ${\cal V}={\cal N}$ and the dimension of ${\cal V}$ is higher than the 
expected one.

\begin{conjecture}
1) Is it true that for negative indices of rigidity the dimension of the 
variety of $(p+1)$-tuples 
with non-trivial centralizers is always smaller than the expected dimension 
of the variety of all $(p+1)$-tuples (of matrices $M_j$ or $A_j$) ? 

2) Is it true that for negative indices of rigidity if the Jordan normal 
forms $J^n_1$, 
$\ldots$, $J^n_{p+1}$ satisfy the conditions of Theorem~\ref{generic}, then 
the Deligne-Simpson problem is weakly solvable for any eigenvalues ?
\end{conjecture}

{\Large\bf Appendix. Proof of Theorem~\ref{Gabber} (by Ofer Gabber)}\\ 

$1^0$. We use arguments related to the ones from \cite{Ka}. Suppose we are 
given the conjugacy classes $C_i\subset GL(n,{\bf C})$, $1\leq i\leq p+1$, 
and we are interested in solutions of 

\begin{equation}\label{star}
M_1\ldots M_{p+1}={\rm id}~,~M_i\in C_i
\end{equation}

We say that a solution $M=(M_1,\ldots ,M_{p+1})$ is {\em rigid} if every 
solution $M'$ in some neighbourhood of $M$ is $GL(n,{\bf C})$-conjugate 
to $M$. Here ``neighbourhood'' can be taken in the classical or in the Zariski 
topology. 

$2^0$. Consider distinct points $a_1$, $\ldots$, 
$a_{p+1}\in {\bf P}_{{\bf C}}^1$ and set $U={\bf P}_{{\bf C}}^1\backslash 
\{ a_1,\ldots ,a_{p+1}\}$. Choose a base point $x_0\in U$ and a standard set 
of generators $\gamma _i\in \pi _1(U,x_0)$ where $\gamma _i$ is freely 
homotopic to a positive loop around $a_i$, $\gamma _1\ldots \gamma _{p+1}=1$ 
(using $\pi _1$ conventions as in Deligne LNM 163).

Then a solution of (\ref{star}) determines a local system $L$ on $U$, 
$L_{x_0}\simeq {\bf C}^n$; the local monodromies are given by the matrices 
$M_i$.
  
$3^0$. Recall that if $f:X\rightarrow Y$ is an algebraic map of irreducible 
algebraic varieties, then every irreducible component of a fibre of $f$ has 
dimension $\geq$dim$(X)-$dim$(Y)$. 

Suppose we are given a rigid solution of (\ref{star}). In particular, 
if $\delta _i$ 
is the value of the determinant on $C_i$, then $\prod \delta _i=1$, so we 
have the product morphism 

\[ f:C_1\times \ldots \times C_{p+1}\rightarrow SL(n,{\bf C})\]
and by assumption the $GL(n,{\bf C})$-orbit of $(M_1,\ldots ,M_{p+1})$ is 
dense in an irreducible component of $f^{-1}($id$)$. The above orbit is also 
an $SL(n,{\bf C})$-orbit, so it is of dimension $\leq n^2-1$. 

$4^0$. Hence, 

\[ \sum _{i=1}^{p+1}d_i\leq 2(n^2-1)~.\] 
Denote by $j$ the inclusion of $U$ in ${\bf P}_{{\bf C}}^1$ and by 
${\cal Z}(M_i)$ the 
space of matrices commuting with $M_i$. Then $d_i=n^2-$dim${\cal Z}(M_i)$ and 
by the Euler-Poincar\'e formula (cf. \cite{Ka} p. 16) the above inequality 
is equivalent to 

\[ \chi ({\bf P}_{{\bf C}}^1,{\rm j}_*\underline{{\rm End}}(L))\geq 2~.\]

Now if $F$ is a rank $n$ irreducible local system with local monodromies in 
the prescribed conjugacy classes, then by the Euler-Poincar\'e formula

\[ \chi ({\bf P}_{{\bf C}}^1,{\rm j}_*\underline{{\rm End}}(L))=
\chi ({\bf P}_{{\bf C}}^1,{\rm j}_*\underline{{\rm Hom}}(L,F))\geq 2~,\]
so one of the the two cohomology groups 
$H^0({\bf P}_{{\bf C}}^1,{\rm j}_*\underline{{\rm Hom}}(L,F))\cong 
{\rm Hom} _U(L,F)$ or  

\[ H^2({\bf P}_{{\bf C}}^1,{\rm j}_*\underline{{\rm Hom}}(L,F))\cong 
H_c^2(U,\underline{{\rm Hom}}(L,F))\cong {\rm Hom} _U(F,L)^{{\sf v}}\]
is non-zero, which implies (as $F$ is irreducible) that $F\simeq L$ 
(cp. \cite{Ka}, Theorem 1.1.2). Hence, if $L$ is reducible, then $F$ does not 
exist.\hspace{1cm}$\Box$

Author's address: Universit\'e de Nice, Laboratoire de Math\'ematiques, 
Parc Valrose, 06108 Nice Cedex 2, France,

tel.: (0033) 4 92 07 62 67

fax : (0033) 4 93 51 79 74

e-mail kostov@math.unice.fr
\end{document}